\documentclass[letterpaper,12pt]{amsart}

\usepackage{amsfonts,amsmath,amsthm,amssymb, graphicx,bbm,nicefrac}

\newtheorem{thm}{Theorem}
\newtheorem{lem}[thm]{Lemma}

\newtheorem{cor}[thm]{Corollary}

\newtheorem{conj}[thm]{Conjecture}

\newcommand{\Phalf}[1]{P_{1/2}(#1)}
\newcommand{\D}[2]{\mathcal{D}_{#1, #2}}
\newcommand{\q}{\mathbf{q}}
\newcommand{\qd}[2]{\mathbf{q}^{\mathbf{d}(#1,#2)}}
\newcommand{\qdinv}[2]{\mathbf{q}^{-\mathbf{d}(#1,#2)}}
\newcommand{\GamCv}[2]{\Gamma_{#1}(#2)}
\newcommand{\GamCjv}[3]{\Gamma_{#1}^{#2}(#3)}
\newcommand{\Pnd}[2]{\mathcal{P}_{#1}^{#2}}
\newcommand{\hCq}[2]{h(#1,#2)}

\def\Pr{\mathbb{P}}

\begin{document}

\title[Threshold Pebbling]{Thresholds for Pebbling on Grids}

\date{}

\author[N.~Bushaw, N.~Kettle]{Neal Bushaw \and
	Nathan Kettle}
\address{(NB): Virginia Commonwealth University, Richmond VA} \email{nobushaw@vcu.edu}
\address{(NK): Instituto de Matem\'atica Pura e Aplicada, Rio de Janeiro BR}

\begin{abstract}
Given a connected graph $G$ and a configuration of $t$ pebbles on the vertices of G, a $q$-pebbling step consists of removing $q$ pebbles from a vertex, and adding a single pebble to one of its neighbours. Given a vector $q=(q_1,\ldots,q_d)$, $q$-pebbling consists of allowing $q_i$-pebbling in coordinate $i$. A distribution of pebbles is called solvable if it is possible to transfer at least one pebble to any specified vertex of $G$ via a finite sequence of pebbling steps.

In this paper, we determine the weak threshold for $\q$-pebbling on the sequence of grids $[n]^d$ for fixed $d$ and $\q$, as $n\to\infty$. Further, we determine the strong threshold for $q$-pebbling on the sequence of paths of increasing length.  A fundamental tool in these proofs is a new notion of `centralness' and a sufficient condition for solvability based on the well used pebbling weight functions; we believe this to be the first result of its kind, and may be of independent interest.

These theorems improve recent results of Czygrinow and Hurlbert, and Godbole, Jablonski, Salzman, and Wierman. They are the generalizations to the random setting of much earlier results of Chung.

In addition, we give a short counterexample showing that the threshold version of a well known conjecture of Graham does not hold.  This uses a result for hypercubes due to Czygrinow and Wagner.
\end{abstract}

%\subjclass[2000]{Primary 05C35, 05C38}

\keywords{pebbling, grids, paths, threshold}

\maketitle

\section{History and Fundamentals}
In this paper, we discuss the so-called 'pebbling game', first introduced by Fan Chung Graham in 1989 \cite{chung1989pebbling}. In this paper, she credits this game to a suggestion of Jeff Lagarias and Saks. This suggestion was based in application - it arose in relation to a problem of Erd\H{o}s and Lemke (see, e.g., \cite{EK}. However, this game is another in a long tradition of marble \& board games dating back thousands of years (see, e.g., \cite{Fraenkel}).

Given a graph $G$ of order $n$, we use $D$ to denote a \emph{configuration} of $t$ unlabeled pebbles on the vertices of $G$; formally, $D$ is a multiset of $V(G)$. The {\emph{pebbling move}}, then, consists of taking precisely two pebbles from any vertex $v$, and placing a single pebble onto any neighbor of $v$; the other is `lost in transit'.

Given $v\in V(G)$, We call a configuration \emph{$v$-solvable} if it is possible to move a pebble to the vertex $v$ via a sequence of pebbling moves. A configuration is \emph{solvable} if it is $v$-solvable for each $v\in V(G)$.

The \emph{pebbling number} of a graph $\mathcal{G}$, $\Pi(\mathcal{G})$, is the minimum number of pebbles such that each initial configuration is solvable.

As mentioned in the first paragraph, pebbling problems date back several decades to a suggestion of Lagarias and Saks as a method for solving a number theoretic problem; this suggestion was carried out by Chung~\cite{chung1989pebbling}, who first mentions this form pebbling in print. In this paper, she proved that the pebbling number of the hypercube, $\Pi(\mathcal{Q}_n)=2^n$. For thorough discussion of $2$-pebbling, we direct the reader to \cite{H}.

It is natural to consider the variation of this problem in which moves in different directions cost a different number of pebbles. Given $\q=(q_1,\ldots,q_d)$, we define $\q$-pebbling on the $d$-dimensional grid $\mathcal{P}_{n_1}~\times~\ldots~\times~\mathcal{P}_{n_d}$ to be the natural generalization of $q$-pebbling in which $q_i$-pebbling steps are used to move a pebble between vertices adjacent on the $i$\textsuperscript{th} co-ordinate. This generalization was first explored by Chung~\cite{chung1989pebbling}, who showed that in the pebbling number for $\q$-pebbling on the $d$-dimensional grid is $\prod_iq_i^{n_i}$.  Throughout, we will consider the grid $\mathcal{P}_{n_1}~\times~\ldots~\times~\mathcal{P}_{n_d}$ as a graph whose typical vertex is denoted by a $d$-tuple $(x_1,\ldots,x_d)$ with each $x_i\in[n_i]=\{1,2,\ldots,n_i\}$.

\section{Results}

In this paper, we are concerned with a randomized version of the pebbling problem. We shall look at the case where our initial distribution of pebbles is chosen uniformly at random from all distributions of $t$ unlabelled pebbles on $V(G)$ (as a sanity check, we point out that there are $\binom{n+t-1}{t}$ such configurations). We denote this probability space of pebble configurations on $V(G)$ by $\D{G}{t}$.

We are interested in the probability that a such a distribution is solvable; that is, how many pebbles are necessary to ensure the initial configuration is solvable with a decent probability? In particular, as is natural and somewhat standard with thresholds, we consider the quantity $$\Phalf{G}\,:=\,\min\left\{k\mid\Pr\left(D\in\D{G}{k}\mbox{ is solvable}\right)\ge\frac{1}{2}\right\}.$$

Bekmetjev, Brightwell, Czygrinow, and Hurlbert~\cite{bekmetjev2003thresholds} showed that for any graph sequence $\mathcal{G}=\left(G_i\right)_{i\ge 1}$, the function $\Phalf{\mathcal{G}_i}$ is a weak threshold for solvability\footnote{It is worth noting that where it introduces no ambiguity, we will often omit details of our graph sequences, and use only the standard Landau-Bachmann asymptotic behavior.}.

In this paper we determine the weak threshold for $\q$-pebbling on $d$-dimensional grids; that is $\left\{\Pnd{n}{d}\right\}_{n\ge1}$. This result appears as Theorem~\ref{main_thm}, and improves a result of Czygrinow and Hurlbert~\cite{czygrinow2006girth} which gave $P_{\frac{1}{2}}(\Pnd{n}{d})=n^d\exp\left(\Theta\left(\log (n^d)\right)^{\frac{1}{d+1}}\right)$. In particular, we prove the following much stronger theorem.

\begin{thm}
\label{main_thm}
For $\q$-pebbling on $\Pnd{n}{d}$, we have that 
$$P_{\frac{1}{2}}\left(\Pnd{n}{d}\right)=n^d\exp\left(\left(\frac{(d+1)!\log n\prod_i\log q_i}{2}\right)^{\frac{1}{d+1}}-\frac{d\log\log n}{(d+1)}+O(1)\right).$$
\end{thm}

We note that, in general, there need not be a strong threshold. As an example, consider the sequence of complete graphs $\mathcal{K}_n$. A configuration here is solvable only when there exists a vertex that initially has at least two pebbles. This event that has probability $1-e^{-c^2}+o(1)$ when $k=c\sqrt{n}$, and thus no sharp threshold exists.

However, for each $q\ge 2$, there is indeed a strong threshold for $q$-pebbling of paths $\mathcal{P}_n$. Again, this improves dramatically on the best known bounds of $P_{\frac{1}{2}}(\mathcal{P}_n)=n\exp\left((\sqrt{\log 2}+o(1))\sqrt{\log n}\right)$, with the upper bound due to Godbole, Jablonski, Salzman, and Wierman~\cite{godbole2004improved}, and the lower bound due to Czygrinow and Hurlbert~\cite{czygrinow2008pebbling}. Indeed, in this case we estimate the $O(1)$ term much better in the previous theorem to obtain the following.

\begin{thm}\label{pathsharp}
 For $q$-pebbling on $\mathcal{P}_n$, we have that
 \begin{equation*}
  \Phalf{\mathcal{P}_n}=n\exp\left(\sqrt{\log q\log n}-\frac{\log\log n}{2}+o(1)\right).
 \end{equation*}
\end{thm}

In the next section we prove only the more general result for grids of arbitrary dimension. We follow this with a discussion of the sharp result for paths (whose proof is nearly identical).

\section{The Weak Threshold for Grids}

In this section, we prove Theorem~\ref{main_thm}. Throughout we shall denote by $N=n^d$ the number of vertices in $\Pnd{n}{d}$, and we shall reserve $k=\lambda N$ for the total number of pebbles placed on the vertices.

We first give a simple estimate of the probability of a given configuration of pebbles.

\begin{lem}
\label{prob_bound}
In the probability space $\D{\Pnd{n}{d}}{\lambda N}$, with $\lambda=o(\sqrt{N})$, for any subset $\{v_1,\ldots v_t\}$ of the vertices with $t=o(\lambda)$, and any quantities $f_1,\ldots f_t$, with sum $\sum_i f_i=m=o(\lambda^2)$, the probability of obtaining $f_i$ pebbles on vertex $v_i$ for each $i$ is
\begin{equation*}
\frac{1}{\lambda^te^{\frac{m}{\lambda}}}(1+o(1)).
\end{equation*}
\end{lem}

\begin{proof}
The probability is equal to
\begin{eqnarray*}
\frac{\binom{N-1-t+k-m}{N-1-t}}{\binom{N-1+k}{N-1}} &=& \frac{(N-1-t+k-m)!k!(N-1)!}{(N-1+k)!(k-m)!(N-1+t)!}\\
&=&(N+k)^{-(t+m)}e^{O\left(\frac{(t+m)^2}{k}\right)}k^me^{O\left(\frac{m^2}{k}\right)}N^te^{O\left(\frac{t^2}{N}\right)}\\
&=&(\lambda+1)^{-t}\left(\frac{\lambda}{\lambda+1}\right)^m(1+o(1))\\
&=&\lambda^{-t}e^{-\frac{m}{\lambda}+O\left(\frac{t}{\lambda}+\frac{m}{\lambda^2}\right)}(1+o(1)). 
\end{eqnarray*}
\end{proof}

Associate the vertex set of $\Pnd{n}{d}$ with $[n]^d$ in the canonical way; we'll also use $\mathbf{e}_i$ for the standard unit vectors representing possible directions for pebbling moves. For two vertices $a=(a_1,\ldots,a_d)$, and $b=(b_1,\ldots,b_d)$ in $\Pnd{n}{d}$, denote by $\mathbf{d}(a,b)$ the vector distance between them; that is $(|a_1-b_1|,\ldots,|a_d-b_d|)$. We further denote by $\qd{a}{b}$ the \emph{pebbling distance} $\prod_iq_i^{|a_i-b_i|}$ - this is equal to the number of pebbles that have to be placed at $b$ in order to move one of them to $a$.

In the following lemma, we give a condition for $v$-solvable configurations in terms of this pebbling distance (often referred to in more general graph case as the pebbling weight).  We believe this result to be the first such sufficient condition for solvability using these weight functions.

\begin{lem}
If a distribution of pebbles $D$ from the probability space $\D{\Pnd{n}{d}}{k}$, is $v$-solvable then
\begin{equation}
\label{fractional}
\sum_{w\in\Pnd{n}{d}}\frac{D(w)}{\qd{w}{v}}\ge 1.
\end{equation}
Also if 
\begin{equation}
\label{greedy}
\sum_{w\in\Pnd{n}{d}}\frac{D(w)}{\qd{w}{v}}>\sum_{w\neq v}\frac{\max_i(q_i)-1}{\qd{w}{v}},
\end{equation}
then $D$ is $v$-solvable.   
\end{lem}

\begin{proof}
\emph{Fractional $\q$-pebbling} is a less restrictive form of pebbling where a move consists of any positive amount $\epsilon$ of pebbles being removed from a vertex $y$, and $\epsilon q_i^{-1}$ pebbles being placed on the vertex $y\pm \mathbf{e}_i$. It is clear that if $D$ is $v$-solvable under $q$-pebbling then it is $v$-solvable under fractional $q$-pebbling. Trivially, a configuration is $v$-solvable under fractional $q$-pebbling exactly when \eqref{fractional} holds. 

If a pebble is moved towards $v$ under $\q$-pebbling then the left hand side of \eqref{greedy} does not change. If \eqref{greedy} holds, then either $v$ has a pebble on it or there exists a vertex $w\neq v$ with at least $\max_i(q_i)$ pebbles on it. In the second case we can move a pebble from $w$ towards $v$, and as our graph and total number of pebbles are finite, we can not repeatedly move pebbles towards $v$ without eventually moving a pebble onto $v$. 
\end{proof}

Using \eqref{greedy}, we easily deduce an equivalent bound in a more convenient form.

\begin{cor}
\label{solvable}
If a distribution of pebbles $D$ from $\D{\Pnd{n}{d}}{k}$ satisfies
\begin{equation*}
\sum_{w\in\Pnd{n}{d}}\frac{D(w)}{\qd{w}{v}}>(\max_i(q_i)-1)\prod_i\left(\frac{q_i+1}{q_i-1}\right),
\end{equation*}
then $D$ is $v$-solvable.
\end{cor}

\begin{proof}
We have that 
\begin{eqnarray*}
\sum_{w}\frac{\max_i(q_i)-1}{\qd{w}{v}}&<&(\max_i\{q_i\}-1)\prod_i\left(1+\frac{2}{q_i}+\frac{2}{q_i^2}+\ldots\right)\\
&=&(\max_i(q_i)-1)\prod_i\left(\frac{q_i+1}{q_i-1}\right).
\end{eqnarray*}
\end{proof}

In the next lemma, we give a general result which will allow us to count the possible number of configurations satisfying the bounds of the preceding lemmas.

\begin{lem}
\label{tetra_region}
Suppose that $(a_i)_{i=1}^r$ is a sequence of positive integers, then the number of solutions in non-negative integers $(x_i)_{i=1}^r$ to
\begin{equation*}
\sum_{i=1}^r\frac{x_i}{a_i}<1,
\end{equation*}
is bounded between
\begin{equation*}
\frac{\left(1-\sum_{i=s+1}^ra_i^{-1}\right)^{r-s}\prod_{i=s+1}^ra_i}{(r-s)!},
\end{equation*}
and
\begin{equation*}
\frac{\left(1+\sum_{i=1}^ra_i^{-1}\right)^r\prod_{i=1}^ra_i}{r!},
\end{equation*}
for any $s< r$.
\end{lem}

\begin{proof}
Write $S$ for the set of solutions $(x_i)_{i=1}^r$. The region $Y$ of $\mathbb{R}^r$ determined by $y_i\ge0$, and $\sum_{i=1}^r (y_i-1)a_i^{-1}<1$ entirely contains the region $X$ equal to the union over $S$ of $[x_1,x_1+1]\times\ldots\times [x_r,x_r+1]$. By comparing volumes, this gives
\begin{equation*}
|S|<\frac{\left(1+\sum_{i=1}^ra_i^{-1}\right)^r\prod_{i=1}^ra_i}{r!}.
\end{equation*}
Similarly the region $Y^\prime$ of $\mathbb{R}^{r-s}$ determined by $y_i\ge 0$ for $s<i\le r$, and $\sum_{i=s+1}^r (y_i+1)a_i^{-1}<1$ is entirely contained in the region $X^\prime$ equal to the union over $(x_i)$ in $S$ with $x_i=0$ for $i\le s$ of $[x_{s+1},x_{s+1}+1]\times\ldots\times[x_r,x_r+1]$. Again by comparing volumes this gives
\begin{equation*}
|S|>\frac{\left(1-\sum_{i=s+1}^ra_i^{-1}\right)^{r-s}\prod_{i=s+1}^ra_i}{(r-s)!}.
\end{equation*}
\end{proof}

Throughout the remainder of the proof, we will call a vertex $v=(v_1,\ldots, v_d)$ $(C,t)$-\emph{central} if it has exactly $t$ coordinates satisfying $\frac{\log C}{\log 2}<v_i<n-\frac{\log C}{\log 2}$. Note that there are $O_d((\log C)^{d-t}n^t)$ $(C,t)$-central vertices in $\Pnd{n}{d}$.  This crucial definition allows us to keep track of those vertices which are far from the boundary in many directions (where we can simplify our arguments since no pebble can reach a boundary vertex).  Appropriate choices of $C$ and $t$ will allow us to ensure that most vertices are, in fact, $(C,t)$-central, and thus utilize these simplified counting techniques.

We further define $\GamCv{C}{v}$ to be the set of vertices $w$ with $\qd{w}{v}\le C$, and $w_i\ge v_i$ if $v_i\le\frac{\log C}{\log 2}$, (and similarly $w_i\le v_i$ if $v_i\ge n-\frac{\log C}{\log 2}$). These definitions are made so that we can estimate quite accurately the size of the $\Lambda_C$ neighborhood of a $(C,t)$-central vertex; the technical statement of this follows.

\begin{cor}
\label{tetra_region_cor}
For $v$ a $(C,t)$-central vertex in $\Pnd{n}{d}$, we have that
\begin{equation*}
 |\GamCv{C}{v}| = \frac{2^t(\log C)^d}{d!\prod_{i=1}^d\log q_i}+O_d\left((\log C)^{d-1}\right).
\end{equation*}

\end{cor}

\begin{proof}
 Let $v_0=(1,\ldots,1)$ - note that $v_0$ is trivially $(C,0)$-central. $\GamCv{C}{v_0}$ consists of all vertices of the form $(x_1+1,\ldots,x_d+1)$ with each $x_i$ non-negative and
 \begin{equation}
  \label{linear_bound}
  \sum_i x_i\log q_i \le \log C.
 \end{equation}
From Lemma~\ref{tetra_region}, we obtain
 \begin{equation*}
  |\GamCv{C}{v_0}|=\frac{(\log C)^d}{d!\prod_{i=1}^d\log q_i}+O\left((\log C)^{d-1}\right).
 \end{equation*}
 For a general $(C,t)$-central vertex $v$, we can see that $\GamCv{C}{v}$ is the union of $2^t$ regions isomorphic to $\GamCv{C}{v_0}$. The intersection of any two of these regions must be contained $\GamCv{C}{v}\cap\{u\mid u_i=v_i\}$ for some $i$. This region is isomorphic to $\GamCv{C}{v^\prime}$ for some $v^\prime\in\Pnd{n}{d-1}$, and hence has size $O_d((\log C)^{d-1})$ by induction on $d$. Putting this together, we see that 
 \begin{equation*}
  |\GamCv{C}{v}| = \frac{2^t(\log C)^d}{d!\prod_{i=1}^d\log q_i}+O_d\left((\log C)^{d-1}\right).
 \end{equation*}
\end{proof}

We next give a bound on the influence of vertices at pebbling distance at least $C$ from $v$.

\begin{lem}
\label{small_sum}
 For $v\in\Pnd{n}{d}$, we have that
 \begin{equation*}
  \sum_{w:\qd{w}{v}>C}\qdinv{w}{v}=O_d\left(\frac{(\log C)^{d-1}}{C}\right).
 \end{equation*}

\end{lem}
\begin{proof}
 The sum is maximized for a $(C,d)$-central vertex $v$; we can bound the sum by
 \begin{equation*}
  \sum_{w:\qd{w}{v}>C}\qdinv{w}{v}\le\sum_{j=0}\frac{|\GamCv{2^{j+1}C}{v}\setminus\GamCv{2^jC}{v}|}{2^jC}.
 \end{equation*}
 
 Noting that  
 \begin{eqnarray*}
  |\GamCv{2^{j+1}C}{v}\setminus\GamCv{2^jC}{v}|&\le& \frac{2^d((\log 2^{j+1}C)^d-(\log 2^jC)^d)}{d!\prod_{i=1}^d\log q_i}+O_d((\log 2^jC)^{d-1})\\
  &=& O_d\left((\log 2^jC)^{d-1}\right),
 \end{eqnarray*}
 the claim follows without effort.
\end{proof}

Combining the preceding, we attain the following estimate.

\begin{lem}
\label{product}
For a $(C,t)$-central vertex $v$, we have that
\begin{equation*}
\prod_{w\in\GamCv{C}{v}}\qd{v}{w}=\exp\left(\frac{d2^{t}(\log C)^{d+1}}{(d+1)!\prod_{i=1}^d\log q_i}+O_d((\log C)^d)\right).
\end{equation*}
\end{lem}

\begin{proof}
When $d=1$, we have that $\GamCv{C}{v}$ contains $t+1$ vertices at every distance from $v$, up to $\left\lfloor\frac{\log C}{\log q_1}\right\rfloor$, and so
\begin{eqnarray*} 
\prod_{w\in\GamCv{C}{v}}\qd{v}{w}&=&\exp\left(\frac{(t+1)\log q_1}{2}\left(\left\lfloor\frac{\log C}{\log q_1}\right\rfloor+\left\lfloor\frac{\log C}{\log q_1}\right\rfloor^2\right)\right)\\
&=&\exp\left(\frac{(t+1)(\log C)^2}{2\log q_1}+O(\log C)\right),
\end{eqnarray*}
which is equivalent to the stated equality.

For $d>1$, we can partition $\GamCv{C}{v}$ according to the value of its $d$\textsuperscript{th} coordinate. Let $\GamCjv{C}{j}{v}=\GamCv{C}{v}\cap\{x\mid x_d=v_d+j\},$ and set $\eta=1$ when $(\frac{\log C}{\log 2})<v_d<n-\frac{\log C}{\log 2})$, but $\eta=0$ otherwise. Notice that as the vertices in $\GamCjv{C}{j}{v}$ are those that have $\prod_{i=1}^{d-1}q_i^{|w_i-v_i|}\le\frac{C}{q_d^j}$, we have 
\begin{equation*}
\prod_{w\in\GamCjv{C}{j}{v}}\qd{v}{w}=\exp\left(|j|\log q_d\left|\GamCjv{C}{j}{v}\right|+\frac{(d-1)2^{t-\eta}(\log C-|j|\log q_d)^{d}}{d!\prod_i^{d-1}\log q_i}+O_d((\log C)^{d-1})\right),
\end{equation*}
and, by Corollary~\ref{tetra_region_cor},
\begin{equation*}
|\GamCjv{C}{j}{v}|=\frac{2^{t-\eta}(\log C-|j|\log q_d)^{d-1}}{(d-1)!\prod_{i=1}^{d-1}\log q_i}+O_d\left((\log C)^{d-2}\right).
\end{equation*}
We therefore have
\begin{equation*}
\prod_{w\in\GamCv{C}{v}}\qd{v}{w}=\exp\left(\frac{2^{t-\eta}}{d!\prod_{i=1}^{d-1}\log q_i}\sum_{j=-\lfloor\frac{\log C}{\log q_d}\rfloor}^{\lfloor\frac{\log C}{\log q_d}\rfloor}\beta_j+O\left((\log C)^d\right)\right),
\end{equation*}
where
\begin{eqnarray*}
 \beta_j & = & d|j|\log q_d(\log C-|j|\log q_d)^{d-1}+(d-1)(\log C-|j|\log q_d)^d\\
	 & = & \sum_{i=0}^{d}d\binom{d-1}{i-1}(-1)^{i-1}|j|^i(\log q_d)^i(\log C)^{d-i}+(d-1)\binom{d}{i}(-1)^{i}|j|^i(\log q_d)^i(\log C)^{d-i}\\
	 & = & \sum_{i=0}^d(d-i-1)\binom{d}{i}(-1)^i|j|^i(\log q_d)^i(\log C)^{d-i}.
\end{eqnarray*}
As
\begin{equation*}
 \sum_{j=-\lfloor t \rfloor}^{\lfloor t\rfloor}|j|^k=\frac{2t^{k+1}}{k+1}+O\left(t^k\right),
\end{equation*}
and
\begin{equation*}
 \sum_{j=0}^{\lfloor t\rfloor}|j|^k=\frac{t^{k+1}}{k+1}+O\left(t^k\right),
\end{equation*}
we have that $\prod_{w\in\GamCv{C}{v}}\qd{v}{w}$ is equal to
\begin{eqnarray*}
\exp\left(\frac{2^{t-\eta}}{d!\prod_{i=1}^{d-1}\log q_i}\sum_{i=0}^d(d-i-1)\binom{d}{i}(-1)^i\frac{2^{\eta}(\log C)^{d+1}}{(i+1)\log q_d}+O_d\left((\log C)^d\right)\right),
\end{eqnarray*}
from which the lemma follows.
\end{proof}

At long last, we reach the Lemma which instigated this notion of centrality: for a $(C,t)$ central vertex, we can count the number of configurations on $\GamCv{C}{v}$ in which moving $\ell$ vertices to $v$ is possible.

\begin{lem}
\label{distribution_number}
 Let $v$ be a $(C,t)$-central vertex. The number of distributions $D$, of pebbles on $\GamCv{C}{v}$ such that 
 \begin{equation}
  \sum_{w\in\GamCv{C}{v}}\frac{D(w)}{\qd{v}{w}}<\ell,
 \end{equation}
is
\begin{equation}
 \exp\left(\frac{2^t(\log C)^d}{(d-1)!\prod_{i=1}^d\log q_i}\left(\frac{\log C}{d+1}-\log\log C\right)+O_{d,l}\left((\log C)^d\right)\right).
\end{equation}
\end{lem}

\begin{proof}
Applying Lemma~\ref{tetra_region} with $a_i=\ell\qdinv{v}{w_i}$ for $w_i\in\GamCv{C}{v}$, and $r=|\GamCv{C}{v}|$, we immediately obtain an upper bound of
 \begin{equation}
  \frac{\left(1+\frac{\prod_i\left(\frac{q_i+1}{q_i-1}\right)}{\ell}\right)^{|\GamCv{C}{v}|}\ell^{|\GamCv{C}{v}|}\prod_{w\in\GamCv{C}{v}}\qd{w}{v}}{|\GamCv{C}{v}|!}.
 \end{equation}

 By Corollary~\ref{tetra_region_cor} and Lemma~\ref{product}, this upper bound is equal to
 \begin{equation}
  \exp\left(\frac{d2^t(\log C)^{d+1}}{(d+1)!\prod_{i=1}^d\log q_i}-\frac{2^t(\log C)^d\log|\GamCv{C}{v}|}{d!\prod_{i=1}^d\log q_i}+O_{d,\ell}((\log C)^d)\right).
 \end{equation}
Since $\log|\GamCv{C}{v}|=d\log\log C+O_d(1)$, this is precisely the claimed bound.

 The lower bound is nearly identical; however we first need to note that by Lemma~\ref{small_sum}, there exists a constant $\theta_d$ such that
 \begin{equation}
  \sum_{w\in\GamCv{C}{v}\setminus\GamCv{\theta_d}{v}}\qdinv{w}{v}<\frac{1}{2},
 \end{equation}
and so we can take $s=|\GamCv{\theta_d}{v}|$ in Lemma~\ref{tetra_region}.
\end{proof}

With our machinery in place, we are now ready to prove an upper bound on the threshold for $\q$-pebbling on $\Pnd{n}{d}$.

\begin{proof}[Proof of Theorem~\ref{main_thm} upper bound]

We pebble $\Pnd{n}{d}$ according to $\D{\Pnd{n}{d}}{k}$ with $k=\lambda N$ pebbles, where 
 \begin{equation}
  \lambda=\exp\left(\left(\frac{(d+1)!\log n \prod_i\log q_i}{2}\right)^{\frac{1}{d+1}}-\frac{d\log\log n}{d+1}+\gamma\right),
 \end{equation}
 for some large constant $\gamma$.
 
 By Corollary~\ref{solvable}, we have that letting $\ell=(\max_i(q_i)-1)\prod_i\left(\frac{q_i+1}{q_i-1}\right)$ and $v$ be a $(C,t)$-central vertex,
 \begin{equation*}
  \Pr(\Pnd{n}{d}\mbox{ not $v$-solvable})\le\Pr\left(\sum_{w\in\GamCv{C}{v}}\qdinv{w}{v}<\ell\right).
 \end{equation*}

Set $C=\lambda(\log\lambda)^d$. Every distribution of pebbles on $\GamCv{C}{v}$, with $\sum_{w\in\GamCv{C}{v}}\qdinv{w}{v}<\ell$, has at most $lC=o(\lambda^2)$ pebbles. Therefore by Lemma~\ref{prob_bound}, each such distribution has a probability of occuring which is at most $\lambda^{-|\GamCv{C}{v}|}$. Summing over all such distributions and recalling Lemma~\ref{distribution_number}, we see that

 \begin{eqnarray*}
  \Pr\left(\sum_{w\in\GamCv{C}{v}}\qdinv{w}{v}<\ell\right)&\le&\exp\bigg(\frac{2^t(\log C)^d}{(d-1)!\prod_{i=1}^d\log q_i}\left(\frac{\log C}{d+1}-\log\log C-\frac{\log\lambda}{d}\right)\\&&\qquad\qquad+O_{d,\ell}\left((\log C)^d\right)\bigg)\\
  &=&\exp\left(-\frac{2^t(\log C)^{d+1}}{(d+1)!\prod_{i=1}^d\log q_i}+O_{d,\ell}\left((\log C)^d\right)\right).
 \end{eqnarray*}
 By our choice of $C$, we have

\begin{equation*}
  (\log C)^{d+1}=\frac{(d+1)!\log n\prod_i\log q_i}{2}+\gamma(\log C)^d+O_{d,l}(\log C)^d,
\end{equation*}
and so in combination with the preceding inequality, for $\gamma$ sufficiently large, we see that
\begin{equation*}
 \Pr\left(\sum_{w\in\GamCv{C}{v}}\qdinv{w}{v}<l\right)\le\exp\left(-2^{t-1}\log n-(\log n)^{\frac{d}{d+1}}\right).
\end{equation*}

Taking a union bound over all $(C,t)$-central vertices we then conclude
\begin{equation}
 \Pr(\Pnd{n}{d}\mbox{ is not solvable})\le\sum_{t=0}^d2^{d-t}\binom{d}{t}n^t\left(\frac{\log C}{\log 2}\right)^{d-t}\exp\left(-2^{t-1}\log n-(\log n)^{\frac{d}{d+1}}\right),
\end{equation}
which is $o(1)$ as $2^{t-1}\ge t$, and $(\log C)^{d-t}\exp\left(-(\log n)^{\frac{d}{d+1}}\right)=o(1)$.
\end{proof}

Now, we proceed to the lower bound. We shall need a one extra lemma; a straightforward application of \ref{prob_bound} gives an upper bound on the size of the largest pile of pebbles on any vertex.

\begin{lem}
\label{max_degree}
In $\D{\Pnd{n}{d}}{\lambda N}$, with $\log N \ll \lambda \ll \sqrt{N}$, the maximum number of pebbles on any vertex is at most $(1+o(1))\lambda\log N$ with high probability.
\end{lem}

\begin{proof}
By Lemma~\ref{prob_bound} the probability a vertex has more than $(1+\epsilon)\lambda\log N$ pebbles is at most 

\begin{equation*}
(1+o(1))\frac{1}{\lambda e^{(1+\epsilon)\log N}}\left(1+e^{-\frac{1}{\lambda}}+e^{-\frac{2}{\lambda}}+\ldots\right)\le (1+o(1))N^{-(1+\epsilon)}.
\end{equation*}

A union bound over all vertices in $\Pnd{n}{d}$ now gives the result. 
\end{proof}

We now proceed to the proof of the lower bound in Theorem~\ref{main_thm}.

\begin{proof}[Proof of Theorem~\ref{main_thm} lower bound]
We again pebble $\Pnd{n}{d}$ with $k=\lambda N$ pebbles according to $\D{\Pnd{n}{d}}{k}$, with 
 \begin{equation}
  \lambda=\exp\left(\left(\frac{(d+1)!\log n \prod_i\log q_i}{2}\right)^{\frac{1}{d+1}}-\frac{d\log\log n}{d+1}-\gamma\right),
 \end{equation}
 for some large constant $\gamma$.

Let $S$ be the set of vertices in $\Pnd{n}{d}$ with every coordinate equal to $n$, except the $d$th coordinate, which is between $\frac{\log C^\prime}{\log 2}$ and $n-\frac{\log C^\prime}{\log 2}$; the vertices in $S$ are $(C^\prime,1)$-central, with $C^\prime=\lambda(\log N)^2$. We will show that with high probability $S$ contains many vertices satisfying
\begin{equation}
\label{hard_to_pebble}
 \sum_{w\in\GamCv{C}{v}}D(w)\qdinv{w}{v}<\frac{1}{2},
\end{equation}
when $C=\lambda(\log n)^{1-1/d}\log\log n$.

Let $p$ be the probability a vertex $v\in S$ satisfies \eqref{hard_to_pebble}; note that this is the same for each $v\in S$. As every pebbling of $\GamCv{C}{v}$ satisfying \eqref{hard_to_pebble} contains at most $\frac{C}{2}=o(\lambda^2)$ pebbles, we can use Lemma~\ref{prob_bound} to show that any such pebbling has probability at least
\begin{equation*}
 \exp\left(-\frac{2(\log C)^d\log\lambda}{d!\prod_{i=1}^d\log q_i}-\frac{C}{2\lambda}+O_d(\log\lambda(\log C)^{d-1})\right).
\end{equation*}

Since $\frac{C}{\lambda}=\log n^{1-1/d}\log\log n=o((\log n)^{1-1/(d+1)})$ and $\log\lambda\sim\log C$, a union bound over all pebblings of $\GamCv{C}{v}$ thus gives us a lower bound of
\begin{equation*}
 p\ge \exp\left(\frac{2(\log C)^d}{(d-1)!\prod_{i=1}^d\log q_i}\left(\frac{\log C}{d+1}-\log\log C-\frac{\log\lambda}{d}+O_d(1)\right)\right).
\end{equation*}

Substituting in $C=\lambda(\log n)^{1-1/d}\log\log n$, we see that
\begin{equation*}
 p\ge\exp\left(-\frac{2\left((\log\lambda)^{d+1}+d(d+1)\log\log\lambda(\log\lambda)^d\right)}{(d+1)!\prod_{i=1}^d\log q_i}+O_d\left((\log\lambda)^d\right)\right).
\end{equation*}
Taking $\gamma$ sufficiently large, we then get a bound of 
\begin{equation*}
 p\ge\exp\left(-\log n+(\log n)^{\frac{d}{d+1}}\right).
\end{equation*}

We define random variables $X_v$, each of which is the indicator of the event that $v$ satisfies \eqref{hard_to_pebble}, and we let $X_S$ be the sum over vertices in $S$ of $X_v$. We have that $\mathbb{E}(X_S)=|S|p\ge\exp\left((\log n)^{d/(d+1)}\right)$. For two vertices $w,v\in S$, with $d$\textsuperscript{th} coordinate differing by at least $\frac{\log C}{\log 2}$, $\GamCv{C}{v}$, and $\GamCv{C}{w}$ are distinct. By Lemma~\ref{prob_bound}, we have for any pebblings $P_v$ of $\GamCv{C}{v}$ and $P_w$ of $\GamCv{C}{w}$ satisfying \eqref{hard_to_pebble},
\begin{equation*}
 \Pr(P_w\cap P_v)=\Pr(P_w)\Pr(P_v)(1+o(1)).
\end{equation*}

Summing over all such pairs of pebblings, we see that
\begin{equation*}
 \Pr(X_v\cap X_w)=\Pr(X_v)\Pr(X_u)(1+o(1))=p^2(1+o(1)).
\end{equation*}

By Chebyshev's inequality (see, e.g., \cite{MGT}), the probability that $X_S$ is at most $\frac{\mathbb{E}(X_S)}{2}$ is at most
\begin{eqnarray}
 \frac{4\mbox{Var}(X_S)}{\mathbb{E}(X_S)^2}&=&\frac{4\sum_{v,w\in S}\Pr(X_v\cap X_w)}{p^2|S|^2}-4\nonumber\\
 &=&\frac{4\left(\sum_{(v,w)\in S^\prime}\Pr(X_v\cap X_w)+\sum_{(v,w)\in S^2\setminus S^\prime}\Pr(X_v\cap X_w)\right)}{p^2|S|^2}-4\nonumber,
\end{eqnarray}
where $S^\prime$ consists of those pairs of vertices from $S$ whose $d$\textsuperscript{th} coordinates differ by at most $\frac{\log C}{\log 2}$. Bounding the probabilities in these sums by $p$ and $p^2(1+o(1))$ respectively we get that the probability that $X_S$ is at most $\frac{\mathbb{E}(X_S)}{2}$ is at most
\begin{eqnarray}
 \frac{4p|S^\prime|+4p^2|S^2\setminus S^\prime|(1+o(1))}{p^2|S|^2}-4&=& \frac{4(1-p)|S^\prime|}{p|S^2|}+o(1)\nonumber\\
 &\le&\frac{4n\frac{2\log C}{\log 2}}{n\exp\left((\log n )^{\frac{d}{d+1}}\right)}+o(1)\nonumber\\
 &=&o(1). 
\end{eqnarray}
Therefore with high probability there are at least $\Theta\left(\exp\left((\log n)^{\frac{d}{d+1}}\right)\right)$ vertices in $S$ satisfying \eqref{hard_to_pebble}.

Let $X^\prime_v$ be the event that 
\begin{equation}
\label{harder_to_pebble}
 \{v\in\GamCv{C^\prime}{v}\setminus\GamCv{C}{v}:D(v)>\lambda\log\log n\}=\emptyset.
\end{equation}

As $|\GamCv{C^\prime}{v}|\lambda\log\log n=o(\lambda^2)$, we can apply Lemma~\ref{prob_bound} to get that for any two distributions of pebbles $P_v$, and $P^\prime_v$ of $\GamCv{C}{v}$, and $\GamCv{C^\prime}{v}\setminus\GamCv{C}{v}$ satisfying \eqref{hard_to_pebble} and \eqref{harder_to_pebble} respectively,
\begin{equation}
 \Pr(P_v\cap P^\prime_v)=\Pr(P_v)\Pr(P^\prime_v)(1+o(1)).
\end{equation}
Summing over all such pairs of distibutions of pebbles, gives
\begin{equation}
 \Pr(X_v\cap X^\prime_v)=\Pr(X_v)\Pr(X^\prime_v)(1+o(1)),
\end{equation}
or equivalently
\begin{equation}
 \Pr(X^\prime_v\mid X_v)=\Pr(X^\prime_v)(1+o(1)).
\end{equation}
We can bound $\Pr(X^\prime_v)$ by summing $\Pr(D(w)>\lambda\log\log n)$ over all vertices in $\GamCv{C^\prime}{v}\setminus\GamCv{C}{v}$, to get
\begin{equation}
 \Pr(X^\prime_v)=O\left((\log n)^{\frac{d}{d+1}}(\log n)^{-1}\right)=o(1).
\end{equation}

By Lemma~\ref{max_degree} we have that with high probability every vertex in $\Pnd{n}{d}$ has at most $2\lambda\log N$ pebbles on it, we also have that with high probability there exists a vertex $v\in S$ satisfying \eqref{hard_to_pebble}, by the previous remark, we have that with high probability $v$ also satisfies \eqref{harder_to_pebble}, and so with high probability there exists $v$ satisfying

\begin{align*}
 \sum_{w}D(w)\qdinv{w}{v}&=\sum_{w\in\GamCv{C}{v}}D(w)\qdinv{w}{v}+\sum_{w\in\GamCv{C^\prime}{v}\setminus\GamCv{C}{v}}D(w)\qdinv{w}{v}\\
 &\qquad+\sum_{w\not\in\GamCv{C^\prime}{v}}D(w)\qdinv{w}{v}\\
 &\le \frac{1}{2}+\lambda\log\log n O_d\left(\frac{(\log C)^d}{C}\right)+2\lambda\log N O_d\left(\frac{(\log C^\prime)^d}{C^\prime}\right)\\
 &= o(1),
\end{align*}

and so $\Pnd{n}{d}$ is not $v$-solvable.
\end{proof}

As stated in the introduction, for the path $\mathcal{P}_n$, we are able to prove the much stronger Theorem \ref{pathsharp}.  The proof is identical to that of Theorem \ref{main_thm}, with the small modification explained below.  As this was also proven independently by Moews in \cite{Moews}, we omit the (redundant) full proof here.
 
 For a distribution $D$ of pebbles, a vertex $i\in [n]$ is solvable if and only if
 \begin{equation*}
  \sum_{j>i}\frac{D(j)}{q^{j-i}}\geq 1,
 \end{equation*}
or
\begin{equation*}
 \sum_{j<i}\frac{D(j)}{q^{i-j}}\geq 1.
\end{equation*}

Therefore, if $i$ is $(C,1)$-central, then the number of distributions of pebbles on $\GamCv{C}{i}$ such that $\mathcal{P}_n$ is $i$-solvable from just those pebbles is equal to $\hCq{\lfloor{\log C/\log q}\rfloor}{q}^2$, where $\hCq{t}{q}$ is the number of partitions of $q^t$ into powers of $q$. The asymptotics of the function $h$ have been well-studied, in particular Mahler~\cite{mahler1940special} showed that,

\begin{align*}
 h(t,q)=\exp\large(&\frac{(t-1)\log q}{2}-(t-1)\log(t-1)\log\log q\\&+\left(\frac{\log q}{2}+1+\log\log q\right)(t-1)+O(\log(t)^2)\large).
\end{align*}

This improved result on the number of such configurations is enough to tame the error term in the proof of Theorem \ref{main_thm}.

\section{The Product Conjecture}
For a graph $G$, we denote by $\pi(G)$ the smallest size of a solvable configuration of pebbles in $2$-pebbling.  Graham's Conjecture is a well known conjecture stating simply that for any two graphs $G_1$ and $G_2$, the Cartesian product $G_1\square G_2$ ought to have a solvable configuration using at most $\pi(G_1)\cdot\pi(G_2)$ pebbles. This well known conjecture has inspired a great deal of work since it was first mentioned in \cite{chung1989pebbling}.

It is easy to translate Graham's Conjecture into our probabilistic context; such conjectures appear in, e.g., \cite{czygrinow2006girth} where they prove that this does hold in the case of paths.  The translated conjecture follows.

\begin{conj}\label{threshgraham} There is a universal constant $C$ such that for any two graphs $G_1$ and $G_2$, we have 
\begin{equation}\label{threshgrahameq} \Phalf{G_1\times G_2}\le C\Phalf{G_1}\cdot\Phalf{G_2}.
\end{equation}
\end{conj}

Consider the sequence $(Q_1,Q_2,\ldots)$, where $Q_i$ is the $i$-dimensional hypercube. The following was proven in \cite{CzWag}.

\begin{thm}\label{cubes}
For every $\varepsilon>0$, $\Phalf{Q_n}=O\left(\frac{2^n}{\left(\log n\right)^{1-\varepsilon}}\right)$.
\end{thm}

Using this result, we give a short counterexample to Conjecture \ref{threshgraham}.

\begin{proof}[Counterexample to Conjecture \ref{threshgraham}]
Consider the sequence $(Q_1,Q_2,\ldots)$, where $Q_i$ is the $i$-dimensional hypercube. As $\Phalf{Q_n}=o(2^n)$, letting $C$ be the universal constant in the conjecture, there is some $n_0$ such that $\Phalf{Q_{n_0}}\le\frac{2^{n_0}}{2C}$.

The reason that hypercubes are convenient for us in this counterexample is that for any $a,b\in\mathbb{N}$, $Q_{a+b}=Q_{x}\times Q_{y}$. Thus if Conjecture \ref{threshgraham} held, by repeatedly applying \eqref{threshgrahameq}, we get that for any $s>0$
\[\Phalf{Q_{2^sn_0}}\le C^{2^s-1}\Phalf{Q_{n_0}}^{2^s}\le\frac{s^{2^s n_0}}{2^{2^s}}=2^{\left(1-\frac{1}{n_0}\right)2^sn_0}.\]
Since $n_0$ is fixed, this is a contradiction to the fact that for every $\varepsilon>0$, $\Phalf{Q_n}=\Omega(2^{n(1-\varepsilon)})$, and thus Conjecture \ref{threshgraham} cannot hold.
\end{proof}

\section*{Acknowledgement}
The authors wish to thank Glenn Hurlbert for his many helpful comments on this manuscript.

\bibliography{Pebblingbib}
\bibliographystyle{amsplain}

\end{document}